\documentclass{article}
\usepackage{amssymb}

\usepackage{amsmath}
\usepackage{chicago}

\newtheorem{theorem}{Theorem}[section]
\newtheorem{proposition}[theorem]{Proposition}

\makeatletter
\@addtoreset{equation}{section}
\makeatother

\input{tcilatex}

\begin{document}

\title{On the Chernoff bound for efficiency of \\
quantum hypothesis testing}
\author{Vladislav Kargin\thanks{%
Cornerstone Research, 599 Lexington Avenue, New York, NY 10022, USA;
slava@bu.edu} }
\date{}
\maketitle

\begin{abstract}
The paper estimates the Chernoff rate for the efficiency of quantum
hypothesis testing. For both joint and separable measurements, approximate
bounds for the rate are given if both states are mixed and exact expressions
are derived if at least one of the states is pure. The efficiency of tests
with separable measurements is found to be close to the efficiency of tests
with joint measurements. The results are illustrated by a test of quantum
entanglement.
\end{abstract}

\section{Introduction}

Mark Kac once called Probability Theory measure theory with a ``soul''
provided by Physics, games of chance, Economics or Geometry.\footnote{%
In a preface to a book about geometric probability by \citeN{santalo76}.} In
a sense then, Quantum Statistics can be called probability theory with a
``subconscious''. The probability distributions, so important for classical
statistics, are no longer the deepest layer of foundations but only an
outward manifestation of geometry in the Hilbert spaces of quantum states.
This foundational change begs for a new look at the classical statistics
results, and this paper contributes by reconsidering the Chernoff-Hoeffding
results about hypothesis testing.

Why quantum statistics? Today, quantum states can be manufactured. For
example, in one method (\citeN{cirac_zoller95}) ions are placed in a trap
created by electrostatic potential and radio-frequency oscillations. The
ions then are cooled by laser emission, and arranged on a line in the trap.
After that, each individual ion can be accessed by laser pulses and their
joint quantum state can be altered according to the researcher's wishes.
This ability to built and manipulate quantum systems is changing our
thinking about computation and information transmission. Suddenly, certain
classic problems -- \ the factorization of large integers, the search in an
unstructured database, secure communication -- are not as difficult as they
always were.

This conceptual change also affects statistics.

For example, how can a quantum state manufacturer check if states have been
generated faithfully? We can anticipate the statistician's answer: Select a
sample of the states and perform a statistical test. But now, besides
designing the test, the statistician must play an additional role, the role
of advisor on how to perform measurements of a sample of quantum states.
Since in quantum mechanics both the measurement and the state determine the
probability distribution of outcomes, the choice of measurement affects the
properties of the statistical test.

Not all measurements are readily available. Sometimes it is possible to
measure sample states jointly, as one large quantum state, and sometimes the
states can only be measured separately and simultaneously. Yet another
possibility is that the states must be measured separately and sequentially.
Finally, sometimes the sample states can only be measured partially, for
example, when each state represents several remote particles that cannot be
measured jointly. Clearly, the efficiency of the optimal test will depend on
which measurements are available. In this paper we will concentrate on joint
and separable independent measurements.

For a single state the problem of quantum hypothesis testing was solved by %
\citeN{holevo76} and \citeN{helstrom76}. In this paper, I consider a
different situation: when the researcher has access to several copies of the
same state but may not be able to measure them jointly.

The problem of testing using a sample of states was also considered in %
\citeN{helstrom76} (see also recent results by \citeN{ogawa_nagaoka00} and %
\citeN{parthasarathy01}). These authors considered only joint measurements
and only the situation when one of the errors may go to zero arbitrarily
slowly. In contrast, I consider a Bayesian version of the problem, in which
the researcher aims to minimize a weighted average of both errors, and I
consider both joint and separable measurements.

When joint measurements are available, the problem of testing using a sample
can be solved by applying the Holevo-Helstrom result to the case of tensor
powers of primary states. In this case, my main results provide useful
bounds on both the expected error when the sample is finite and the rate of
decline in error as the number of sample states grows. The bounds are given
in terms of fidelity distance between quantum hypothesis.

In addition I derive explicit expression for the rate of error decline in
cases if either one of the hypothesis specifies pure state, or the states
specified by the hypotheses commute.

For the separable measurements, I mainly concentrate on the asymptotic case.
If at least one of the hypotheses is pure, then the optimal separable
measurement has the same asymptotic error rate as the joint measurement. If
both hypotheses are mixed, then there is a measurement whose performance is
close to the performance of the joint optimal measurement. Together these
results imply that the loss in efficiency associated with restriction on
available measurements is not large.

This paper contributes only to the theory of quantum hypothesis testing. I
do not touch on another rapidly growing area of research: quantum state
estimation. For recent progress in this area see the review article by %
\citeANP{gill01}.

The rest of the paper is organized as follows. Section 2 gives some basic
information about quantum states and measurements and formulates the problem
of quantum hypothesis testing. Section 3 gives a short summary of the
Chernoff-Hoeffding results about hypothesis testing. Sections 4 and 5
discuss joint and separable measurements, respectively. Section 6 presents
an illustration. And Section 7 concludes.

\section{Quantum Hypothesis Testing}

States of quantum-mechanical objects -- electrons, photons, atoms,
molecules, etc. -- are described by density matrices. A density matrix is a
self-adjoint, non-negative operator of a complex Hilbert space with a trace
of 1. In this paper we will be concerned only with finite-dimensional
Hilbert spaces, so the operator is indeed represented by a finite Hermitian
matrix. A particular case is matrices of rank one. They are projectors on
one-dimensional subspaces and called pure states.

States are not directly observable: they can be measured but the outcome of
a measurement is a random variable. More precisely, measurements are sets of
non-negative operators which are required to add up to the identity
operator. Each operator corresponds to a particular outcome of the
measurement, and the probability of outcome $i$ if the state is $\rho $ and
the measurement is $\{M_{i}\}$ is $tr\left\{ M_{i}\rho \right\} .$ An
important subclass is formed by measurements in which \ the outcomes are
orthogonal projectors: $M_{i}M_{j}=\delta _{ij}M_{i},$where $\delta _{ij}$
is the Dirac delta-function.

We consider the following problem: a researcher is given a sample of $N$
identical quantum states, which are either $\rho _{0}$ or $\rho _{1}$ with
the prior probability $1/2.$ He aims to minimize the average probability of
making an incorrect decision about the state by devising a system of
measurements and a decision rule. Can we safely assume that all measurements
are available to the researcher? No.

While in some situations the researcher can make a joint measurement of the
state that represent the total sample, most often he can do only separate
measurements of each state in the sample. If the measurements are done
independently of each other, then we will call them separable independent
measurements. If the measurements can be done sequentially and the
researcher adjusts the current measurement according to the results obtained
in the previous measurements, then they are separable adaptable measurements.

Sometimes, the researcher is even more restricted. This happens, for
example, if a sample quantum state consists of two spatially remote parts
and the researcher can only measure them separately. Mathematically it means
that the operators of the measurement must be block-diagonal in a certain
basis. This setup may raise interesting statistical issues about
identification of the state properties.

\section{Classical Chernoff-Hoeffding Bounds}

This section reviews results by \citeN{chernoff52}, \citeN{sanov57}, and %
\citeN{hoeffding65} about asymptotic error rates in hypothesis testing. For
details the reader can also consult the book by \citeANP{cover_thomas91}.

Consider two multinomial distributions, $P$ and $Q,$ from one of which a
sample is drawn and provided to a researcher. The researcher's task is to
guess the distribution. The sufficient statistic for this problem is the
empirical distribution of the sample, $X,$ and the decision rule is
specified by two complementary sets, $\mathcal{P}$ and $\mathcal{Q}$, of
probability distributions on outcomes. If $X\in \mathcal{P}$, hypothesis $P$
is accepted; otherwise, $Q$ is accepted. It is assumed that $P\in \mathcal{P}%
,$ and $Q\in \mathcal{Q}.$

If the true probability distribution is $P,$ it is the Sanov theorem that
asymptotically the probability of making an error and accepting $Q$ is 
\begin{equation}
\exp \left[ -ND(\mathcal{Q}||P)\right] ,
\end{equation}%
up to a subexponential factor, where $D(\mathcal{Q}||P)$ is the
Kullback-Leibler distance from $P$ to $\mathcal{Q}$:%
\begin{equation}
D(\mathcal{Q}||P)=\min_{S\in \mathcal{Q}}\sum_{i=1}^{N}s_{i}\ln \frac{p_{i}}{%
s_{i}}.
\end{equation}%
It follows that the average probability of making an error declines
asymptotically with growth in $N$: 
\begin{equation}
R\sim \exp \left[ -N\min \left\{ D(\mathcal{Q}||P),D(\mathcal{P}||Q)\right\} %
\right] .  \label{error_rate1}
\end{equation}%
The maximum of the decline rate over all possible $\mathcal{P}$ and $%
\mathcal{Q}$ is sometimes called the Chernoff information distance between
distributions $P$ and $Q:$%
\begin{equation}
D_{c}(P,Q)=\max_{\mathcal{P},\mathcal{Q}}\min \left\{ D(\mathcal{Q}||P),D(%
\mathcal{P}||Q)\right\}
\end{equation}

Hoeffding proved that the optimal sets $\mathcal{P}$ and $\mathcal{Q}$ can
be determined from the maximum likelihood principle: a distribution $S$
belongs to $\mathcal{P}$ if and only if $D(S||P)\leq D(S||Q).$ Intuitively,
in this case distribution $S$ is more likely to be observed if the true
distribution is $P$ rather than if it is $Q,$ so hypothesis $P$ should be
accepted.

For example, for multinomial distribution we have the following formula for
the probability of error in the optimal test:%
\begin{eqnarray}
R &=&\frac{1}{2}\left\{ \sum_{\substack{ %
p_{1}^{x_{1}}...p_{n}^{x_{n}}<q_{1}^{x_{1}}...q_{n}^{x_{n}}  \\ %
x_{1}+...+x_{n}=N}}p_{1}^{x_{1}}...p_{n}^{x_{n}}+\sum_{\substack{ %
p_{1}^{x_{1}}...p_{n}^{x_{n}}>q_{1}^{x_{1}}...q_{n}^{x_{n}}  \\ %
x_{1}+...+x_{n}=N}}q_{1}^{x_{1}}...q_{n}^{x_{n}}\right\} \\
&=&\frac{1}{2}\left\{ 1-\frac{1}{2}\sum_{x_{1}+...+x_{n}=N}\left|
p_{1}^{x_{1}}...p_{n}^{x_{n}}-q_{1}^{x_{1}}...q_{n}^{x_{n}}\right| \right\}
\\
&=&\frac{1}{2}\left\{ 1-\frac{1}{2}\left\| P_{N}-Q_{N}\right\| _{1}\right\}
\label{classical_error_rate2} \\
&\sim &c\exp \left[ -ND_{c}(P,Q)\right] ,
\end{eqnarray}%
where the sums are taken over possible results of sampling $N$ times from a
multinomial distribution with $n$ outcomes: $x_{1}$ is the number of
outcomes of type $1,$ $x_{2}$ is the number of outcomes of type $2,$ and so
on; $P_{N}$ and $Q_{N}$ are distributions on the sample space induced by the
distributions $P$ and $Q$ on outcomes, and $\left\| \cdot \right\| _{1}$ is
the total variation norm.

It also turns out (see \citeN{cover_thomas91} for derivation) that for the
optimal choice of $\mathcal{P}$ and $\mathcal{Q}$, the probability
distributions $S\in \mathcal{Q}$ and $S^{\prime }\in \mathcal{P}$ that
minimize respectively $D(S||P)$ and $D(S^{\prime }||Q)$ are the same and
given by the following formula:%
\begin{equation}
s_{i}=\frac{p_{i}^{\lambda }q_{i}^{1-\lambda }}{\sum_{j=1}^{N}p_{j}^{\lambda
}q_{j}^{1-\lambda }},  \label{error_sample_distribution}
\end{equation}%
where $\lambda $ is chosen in such a way that $D(S||P)=D(S||Q).$ Knowing
expression (\ref{error_sample_distribution}) we can derive another
expression for the asymptotic probability of error:%
\begin{equation}
\frac{1}{N}\ln R=\min_{0\leq \lambda \leq 1}\log
\sum_{i=1}^{N}p_{i}^{\lambda }q_{i}^{1-\lambda }.
\label{classical_error_rate}
\end{equation}

All these derivations presuppose that $P$ and $Q$ are fixed. In quantum
statistics the researcher has the ability to vary $P$ and $Q$ by choosing
the measurement. How does this change the classical results?

\section{Joint Measurements}

In this section we look at the joint measurements of a sample of quantum
states. The minimal expected error obtained in this case is a lower bound on
the error achievable when the set of measurements is restricted. In
addition, the theory for joint measurements provides a fascinating
counterpart to the classical theory of the Chernoff bounds.

\subsection{Generalities}

Joint measurement of all sample states is by definition a measurement of the
tensor product of the sample states. Thus, in effect we have the problem of
testing two alternative hypotheses about a single -- although huge --
quantum state, the problem that was solved by Holevo and Helstrom (see, for
example, \citeN{holevo01}). In our situation we only need to determine what
additional implications follow from the special structure of the state.

If the hypotheses about the quantum state are given by matrices $\rho _{0}$
and $\rho _{1}$ with prior probability of $1/2,$ then according to the
Holevo-Helstrom result, the optimal measurement is an orthogonal measurement
with $d$ outcomes, where $d$ is the dimension of the Hilbert space and the
outcomes are projectors on the eigenvectors of operator $\rho _{0}-\rho
_{1}. $ The decision is made based on the following rule: If the measurement
outcome corresponds to an eigenvector with a positive eigenvalue$,$ then $%
\rho _{0}$ is chosen; otherwise, $\rho _{1}$ is chosen.

The minimal expected error probability that can be achieved after the
optimal measurement is given by the following formula:%
\begin{equation}
R=\frac{1}{2}\left( 1-\frac{1}{2}\left\| \rho _{0}-\rho _{1}\right\|
_{1}\right) ,  \label{holevo_helstrom}
\end{equation}%
where $\left\| \cdot \right\| _{1}$ denotes the sum of the absolute values
of eigenvalues.

In our case the hypothetical states are tensor powers of the individual
states, $\rho _{0}^{\otimes N}$ and $\rho _{1}^{\otimes N},$ where%
\begin{equation}
\rho _{i}^{\otimes N}\equiv \underset{N}{\underbrace{\rho _{i}\otimes \rho
_{i}\otimes ...\otimes \rho _{i}}.}
\end{equation}%
The number of outcomes in the optimal joint measurement is $d^{N},$ so it
can be enormous for large values of $N.$ The error is 
\begin{equation}
R=\frac{1}{2}\left( 1-\frac{1}{2}\left\| \rho _{0}^{\otimes N}-\rho
_{1}^{\otimes N}\right\| _{1}\right) .
\end{equation}%
Note the similarity with classical expression (\ref{classical_error_rate2}).

What is the asymptotic rate of decline in error? Can we explicitly calculate
the distribution of eigenvalues of $\rho _{0}^{\otimes N}-\rho _{1}^{\otimes
N}?$

Initial moments of this distribution are indeed easy to calculate. Let us
introduce a notation for the moments:%
\begin{equation}
\mu _{n}=:\int_{0}^{1}t^{n}dF(t)=\frac{1}{d^{N}}tr\left( \rho _{0}^{\otimes
N}-\rho _{1}^{\otimes N}\right) ^{n},
\end{equation}%
where $F(t)$ is the discrete probability distribution that puts equal
probability weight on each eigenvalue. Then the following Proposition holds

\begin{proposition}
\begin{equation}
\mu _{n}=\frac{1}{d^{N}}\sum_{\{k_{1},...k_{n}\}}\left( -1\right) ^{\sum
k_{i}}\left( tr\left( \rho _{k_{1}}...\rho _{k_{n}}\right) \right) ^{N},
\end{equation}%
where $\{k_{1},...k_{n}\}$ run over the set of all $n-$sequences of $0$ and $%
1.$
\end{proposition}

\textbf{Proof:} The proposition follows from the non-commutative binomial
expansion of $\left( \rho _{0}^{\otimes N}-\rho _{1}^{\otimes N}\right) ^{n}$
and the fact that $tr(\rho _{k_{1}}^{\otimes N}...\rho _{k_{n}}^{\otimes
N})=\left( tr\left( \rho _{k_{1}}...\rho _{k_{n}}\right) \right) ^{N}.$ QED

The advantage of this formula is that for a fixed $n,$ the calculation is as
easy for large as for small values of $N.$ The difficulty is that the number
of terms in this formula grows exponentially with $n$. Therefore the
standard map from the set of moment sequences to the set of distributions is
impractical. In the next sections we will pursue a different approach to
estimation of $\left\| \rho _{0}^{\otimes N}-\rho _{1}^{\otimes N}\right\|
_{1}.$

\subsection{Special Cases}

To get more insight about the behavior of $\left\| \rho _{0}^{\otimes
N}-\rho _{1}^{\otimes N}\right\| _{1}$, it is useful to consider several
special cases: (1) when both states are pure; and (2) when the density
operators commute. In the first case let $\rho _{0}=\left| \psi
_{0}\right\rangle \left\langle \psi _{0}\right| $ and $\rho _{1}=\left| \psi
_{1}\right\rangle \left\langle \psi _{1}\right| .$\footnote{%
For convenience, we use the Dirac ket-bra notation: the elements of the
Hilbert space are denoted as $\left| \psi \right\rangle ,$ and the linear
functionals on the Hilbert space are denoted as $\left\langle \psi \right| .$
In particular, $\left| \psi _{0}\right\rangle \left\langle \psi _{0}\right| $
is the orthogonal projector on $\left| \psi _{0}\right\rangle .$} Then we
have the following result:

\begin{theorem}
\label{joint_pure}If both states are pure, then the average error
probability is 
\begin{equation}
R=\frac{1}{2}\left( 1-\sqrt{1-\left| \left\langle \psi _{0}\right| \left.
\psi _{1}\right\rangle \right| ^{2N}}\right) .  \label{pure_joint}
\end{equation}%
Asymptotically,%
\begin{equation}
R\sim \frac{1}{4}\left| \left\langle \psi _{0}\right| \left. \psi
_{1}\right\rangle \right| ^{2N}\text{ as }N\rightarrow \infty .
\end{equation}
\end{theorem}

\textbf{Proof}: Because of (\ref{holevo_helstrom}), we need only to prove
that for pure states 
\begin{equation}
\left\| \rho _{0}^{\otimes N}-\rho _{1}^{\otimes N}\right\| _{1}=2\sqrt{%
1-\left| \left\langle \psi _{0}\right| \left. \psi _{1}\right\rangle \right|
^{2N}}.
\end{equation}%
We can write 
\begin{equation}
\left\| \rho _{0}^{\otimes N}-\rho _{1}^{\otimes N}\right\| _{1}=\left\|
\left| \psi _{0}^{\otimes N}\right\rangle \left\langle \psi _{0}^{\otimes
N}\right| -\left| \psi _{1}^{\otimes N}\right\rangle \left\langle \psi
_{1}^{\otimes N}\right| \right\| _{1}.
\end{equation}

Operator $\left| \psi _{0}^{\otimes N}\right\rangle \left\langle \psi
_{0}^{\otimes N}\right| -\left| \psi _{1}^{\otimes N}\right\rangle
\left\langle \psi _{1}^{\otimes N}\right| $ acts nontrivially only in a $2-$%
dimensional space spanned by $\psi _{0}^{\otimes N}$ and $\psi _{1}^{\otimes
N},$ and it is easy to compute the operator eigenvalues in this space. They
are%
\begin{equation}
\pm \sqrt{1-\left| \left\langle \psi _{0}^{\otimes N}\right. \left| \psi
_{1}^{\otimes N}\right\rangle \right| ^{2}}=\pm \sqrt{1-\left| \left\langle
\psi _{0}\right| \left. \psi _{1}\right\rangle \right| ^{2N}}.
\end{equation}%
From this and the fact that all other eigenvalues are zero, the first
equality of the theorem follows. The asymptotic expression follows from the
Taylor series for the square root.

QED.

Now consider the case of commuting $\rho _{0}$ and $\rho _{1}.$ Let the
distributions of eigenvalues be $P$ for $\rho _{0}$ and $Q$ for $\rho _{1}$.

\begin{theorem}
If states commute, the average error probability is asymptotically%
\begin{equation}
R\sim c\exp \left[ -ND_{c}(P,Q)\right]
\end{equation}
\end{theorem}

In other words, the probability of error has exactly the same growth rate as
in the classical case.

\textbf{Proof:} Since the density operators $\rho _{0}$ and $\rho _{1}$
commute, we can choose the basis in which they both are diagonal. In this
basis%
\begin{equation}
\left\| \rho _{0}^{\otimes N}-\rho _{1}^{\otimes N}\right\|
_{1}=%
\sum_{k=(k_{1},...,k_{d})}|p_{1}^{k_{1}}...p_{d}^{k_{d}}-q_{1}^{k_{1}}...q_{d}^{k_{d}}|,
\end{equation}%
where $k$ is a partition of $N,$ and $(p_{1,}...,p_{d})$ and $%
(q_{1},...,q_{d})$ are eigenvalues of $\rho _{0}$ and $\rho _{1},$
respectively. On the right-hand side we have $\left\| P_{N}-Q_{N}\right\|
_{1}$, the distance between two multinomial distributions, $P_{N}$ and $%
Q_{N} $, that arise in repeated trials from distributions $P$ and $Q$.
Therefore, because $R=\frac{1}{2}\left( 1-\frac{1}{2}\left\| \rho
_{0}^{\otimes N}-\rho _{1}^{\otimes N}\right\| _{1}\right) $ in \ the
quantum case and $R=\frac{1}{2}\left( 1-\frac{1}{2}\left\|
P_{N}-Q_{N}\right\| _{1}\right) $ in the classical case, the average errors
and their asymptotic growth rates are the same in the quantum and classical
cases.

QED.

\subsection{Bounds}

Let us now derive some simple bounds on the error probability that follows
from known inequalities. These bounds are useful because they are rather
narrow and easy to compute. The first set of bounds follows from
inequalities between quantum fidelity and probability of error. The second
bound is only applicable to the asymptotic rate of error decline, and it
follows from a quantum analog of Stein's lemma.

Recall that \emph{fidelity} between two states is defined as follows:%
\begin{equation}
F(\rho _{0},\rho _{1})=tr\sqrt{\sqrt{\rho _{0}}\rho _{1}\sqrt{\rho _{0}}},
\end{equation}%
where $\sqrt{X}$ is the unique non-negative definite, Hermitian matrix $Y$
such that $Y^{2}=X.$

\begin{theorem}
\label{joint_bounds_theorem}Probability of error for optimal test with joint
measurement satisfies the following bounds:%
\begin{equation}
\frac{1}{2}\left( 1-\sqrt{1-\left[ F\left( \rho _{0},\rho _{1}\right) \right]
^{2N}}\right) \leq R\leq \frac{1}{2}\left[ F\left( \rho _{0},\rho
_{1}\right) \right] ^{N}.  \label{bounds}
\end{equation}%
Asymptotically,%
\begin{equation}
2\log F\left( \rho _{0},\rho _{1}\right) \lesssim \frac{1}{N}\log R\lesssim
\log F\left( \rho _{0},\rho _{1}\right)
\end{equation}%
If $\rho _{0}$ is pure, $\rho _{0}=\left| \psi _{0}\right\rangle
\left\langle \psi _{0}\right| ,$ the probability of error satisfies a
tighter upper bound:%
\begin{equation}
R\leq \frac{1}{2}\left[ F\left( \rho _{0},\rho _{1}\right) \right] ^{2N}=%
\frac{1}{2}\left\langle \psi _{0}\right| \rho _{1}\left| \psi
_{0}\right\rangle ^{N}.
\end{equation}
\end{theorem}

\textbf{Proof: }The first result follows from inequality (44) in %
\citeN{fuchs_graaf99} applied to the case of the sample of $N$ independent
states, and from the fact that $F(\rho _{0}^{\otimes N},\rho _{1}^{\otimes
N})=\left[ F(\rho _{0},\rho _{1})\right] ^{N}.$ The second one is a
consequence of Exercise 9.21 in \citeN{nielsen_chuang00}. For the reader's
convenience, I include below short proofs of these results.

The Fuchs-Graaf result states that for every pair of quantum states, $\rho
_{0}$ and $\rho _{1},$ it is true that 
\begin{equation}
1-F\left( \rho _{0},\rho _{1}\right) \leq \frac{1}{2}||\rho _{0}-\rho
_{1}||_{1}\leq \sqrt{1-F\left( \rho _{0},\rho _{1}\right) ^{2}}.
\label{fuchs-graaf}
\end{equation}%
These inequalities follow because

(1) $F\left( \rho _{0},\rho _{1}\right) =\min_{P,Q}F(P,Q),$ where
distributions $P$ and $Q$ arise from a measurement of states $\rho _{0}$ and 
$\rho _{1},$ and where $F\left( P,Q\right) =:\sum_{i}\sqrt{p_{i}q_{i}};$

(2) $||\rho _{0}-\rho _{1}||_{1}=\max_{P,Q}||P-Q||_{1},$ where $P$ and $Q$
come from a measurement, and where $||P-Q||_{1}=:\sum_{i}|p_{i}-q_{i}|;$

(3) the corresponding inequality holds for probability distributions%
\begin{equation}
1-F\left( P,Q\right) \leq \frac{1}{2}||P-Q||_{1}\leq \sqrt{1-F\left(
P,Q\right) ^{2}}.
\end{equation}%
Indeed, given (1), (2), and (3), the left-hand inequality in (\ref%
{fuchs-graaf}) follows because 
\begin{eqnarray}
&&1-F\left( \rho _{0},\rho _{1}\right) \underset{(1)}{=}1-F\left( P,Q\right) 
\text{ (for certain }P\text{ and }Q\text{)} \\
&&\underset{(3)}{\leq }\frac{1}{2}||P-Q||_{1}\underset{(2)}{\leq }\frac{1}{2}%
||\rho _{0}-\rho _{1}||_{1}.
\end{eqnarray}%
The right-hand inequality follows similarly.

Result (1) is from \citeN{fuchs_caves95}. (2) is a restatement of the
Holevo-Helstrom result (\ref{holevo_helstrom}). The left-hand inequality in
(3) holds because 
\begin{equation}
\sum_{i}|p_{i}-q_{i}|\geq \sum_{i}\left( \sqrt{p_{i}}-\sqrt{q_{i}}\right)
^{2}\geq 2\left( 1-\sum_{i}\sqrt{p_{i}q_{i}}\right) .
\end{equation}%
The right-hand inequality in (3) holds because 
\begin{eqnarray}
\left( \sum_{i}|p_{i}-q_{i}|\right) ^{2} &=&\left( \sum_{i}\left| \sqrt{p_{i}%
}-\sqrt{q_{i}}\right| \left| \sqrt{p_{i}}+\sqrt{q_{i}}\right| \right) ^{2} \\
&\leq &\sum_{i}\left| \sqrt{p_{i}}-\sqrt{q_{i}}\right| ^{2}\sum_{i}\left| 
\sqrt{p_{i}}+\sqrt{q_{i}}\right| ^{2} \\
&=&4\left( 1-\left( \sum_{i}\sqrt{p_{i}q_{i}}\right) ^{2}\right) .
\end{eqnarray}

To prove the second part of the theorem, we need to prove that if $\rho
_{0}=\left| \psi _{0}\right\rangle \left\langle \psi _{0}\right| ,$ then
there is such a measurement and a decision rule that 
\begin{equation}
R\leq \frac{1}{2}\left\langle \psi _{0}\right| \rho _{1}\left| \psi
_{0}\right\rangle .  \label{nielsen}
\end{equation}%
Take measurement $\{P_{\psi _{0}},I-P_{\psi _{0}}\},$ where $P_{\psi _{0}}$
is the projector on vector $\psi _{0}.$ Then the probabilities of the first
and second outcomes are respectively $1$ and $0$ if the state is $\psi _{0},$
and $\left\langle \psi _{0}\right| \rho _{1}\left| \psi _{0}\right\rangle $
and $1-\left\langle \psi _{0}\right| \rho _{1}\left| \psi _{0}\right\rangle $
if the state is $\rho _{1}.$ Define the decision rule as follows: state $%
\psi _{0}$ is accepted if and only if the first outcome occurs. The expected
error of this rule is $\frac{1}{2}\left\langle \psi _{0}\right| \rho
_{1}\left| \psi _{0}\right\rangle .$

In the case of tensor powers (\ref{nielsen}) becomes 
\begin{equation}
R\leq \frac{1}{2}\left\langle \psi _{0}\right| \rho _{1}\left| \psi
_{0}\right\rangle ^{N}.
\end{equation}

QED.

Note that the lower bound of inequality (\ref{bounds}) binds for pure
states. This can be seen from (\ref{pure_joint}) because for pure states $%
F(\rho _{0},\rho _{1})=\left| \left\langle \psi _{0}\right| \left. \psi
_{1}\right\rangle \right| .$ The upper bound of inequality (\ref{bounds})
binds for certain commuting operators.

Another bound follows from results by \citeN{ogawa_nagaoka00}. Define
quantum relative entropy:%
\begin{equation}
D(\rho _{0}||\rho _{1})=tr\left[ \rho _{0}(\log \rho _{0}-\log \rho _{1})%
\right] .
\end{equation}%
Then the following lower bound on the error rate holds.

\begin{theorem}
\label{ogawa_nagaoka}$\log \max \left\{ D(\rho _{0}||\rho _{1}),D(\rho
_{1}||\rho _{0})\right\} \lesssim \frac{1}{N}\log R$
\end{theorem}

\textbf{Proof: }$R$ is the average of error probabilities of two types. Say, 
$R=\frac{1}{2}R_{1}+\frac{1}{2}R_{2}.$ If both $R_{1}$ and $R_{2}$ satisfy
the inequality, then $R$ also does. \citeANP{ogawa_nagaoka00}, proved that
if one of the error probabilities violates this inequality, then the other
error probability must approach one as the sample size grows, so the
inequality will hold for the average of the error probabilities, $R.$

QED.

For an example of two-dimensional states, the bounds are illustrated in
Figures 1,2,3 and 4. The states in the example are linear combinations of
the Pauli matrices:%
\begin{eqnarray}
\rho _{0} &=&\frac{1}{2}\left( I+a\sigma _{1}\right) , \\
\rho _{1} &=&\frac{1}{2}\left( I+(b\cos \theta )\sigma _{1}+(b\sin \theta
)\sigma _{2}\right) ,
\end{eqnarray}%
where 
\begin{equation}
\sigma _{1}=\left( 
\begin{array}{cc}
0 & 1 \\ 
1 & 0%
\end{array}%
\right) ,\;\sigma _{2}=\left( 
\begin{array}{cc}
0 & i \\ 
-i & 0%
\end{array}%
\right) .
\end{equation}%
Figures 1 and 2 suggests that bound from Theorem \ref{ogawa_nagaoka} is a
good estimate of the error if the sample size is small and underestimates
the error if the sample size is large. Figures 3 and 4 suggest that bound
from Theorem \ref{ogawa_nagaoka} is especially good when the hypotheses are
close to completely mixed state, $\frac{1}{2}I.$

\section{Separable Measurements}

In the previous section we have seen that it is difficult to compute the
optimal joint measurement because of the high dimensionality of the problem
involved. Besides, even if the optimal joint measurement is found, it can
have an enormous number of outcomes, so it is hard to realize it in the
laboratory. In this section we turn our attention to separable independent
measurements. The goal is to show that the efficiency of a separable
measurement with a small number of outcomes is not much smaller than the
efficiency of the optimal joint measurement.

Let us denote the probabilities of the $i-th$ outcome as $p_{i}$ and $q_{i}$
depending on whether the state is $\rho _{0}$ or $\rho _{1}.$ The following
theorem about optimal measurements holds:

\begin{theorem}
All outcomes of an optimal measurement are projectors.
\end{theorem}

\textbf{Proof:} Indeed, if the measurement includes an outcome, $M_{0},$
that is not a projector then it can be represented as a sum of projectors
with non-negative coefficients:%
\begin{equation}
M_{0}=\sum_{i=1}^{n}\alpha _{i}M_{i}.
\end{equation}%
Therefore 
\begin{eqnarray}
p_{0} &=&:tr(M_{0}\rho _{0})=\sum_{i=1}^{n}\alpha _{i}p_{i}, \\
q_{0} &=&:tr(M_{0}\rho _{1})=\sum_{i=1}^{n}\alpha _{i}q_{i}.
\end{eqnarray}%
Since function $x^{\lambda }y^{1-\lambda }$ is concave and homogeneous, we
have the following inequality%
\begin{equation}
p_{0}^{\lambda }q_{0}^{1-\lambda }\geq \sum_{i=1}^{n}\left( \alpha
_{i}p_{i}\right) ^{\lambda }\left( \alpha _{i}q_{i}\right) ^{1-\lambda }.
\end{equation}%
Because of (\ref{classical_error_rate}), this inequality implies that we can
decrease the error by using the set of outcomes $\{M_{i}\}$ instead of $%
M_{0}.$ This contradicts the optimality of the measurement.

QED.

How many outcomes does an optimal measurement have? It turns out that if one
of the states is pure, $\rho _{0}=\left| \psi _{0}\right\rangle \left\langle
\psi _{0}\right| ,$ then only two outcomes is needed - a huge reduction
relative to the $d^{N}$ outcomes needed for the optimal joint measurement.

\begin{theorem}
\label{separable_pure}When one of the states is pure, there is an
asymptotically optimal test with two outcomes in each measurement. The
average error probability of the test satisfies the following bound%
\begin{equation}
R\lesssim \frac{1}{2}\left\langle \psi _{0}\right| \rho _{1}\left| \psi
_{0}\right\rangle ^{N}\text{ as }N\rightarrow \infty .
\end{equation}
\end{theorem}

\textbf{Proof: } Take measurement $\{P_{\psi _{0}},I-P_{\psi _{0}}\},$ where 
$P_{\psi _{0}}$ is the projector on vector $\psi _{0}.$ Then the
probabilities of the first and second outcomes are respectively $1$ and $0$
if the state is $\rho _{0},$ and $\left\langle \psi _{0}\right| \rho
_{1}\left| \psi _{0}\right\rangle $ and $1-\left\langle \psi _{0}\right|
\rho _{1}\left| \psi _{0}\right\rangle $ if the state is $\rho _{1}.$ Define
the decision rule as follows: state $\rho _{0}$ is accepted if and only if
the second outcome never occurred. This rule leads to an error if and only
if the true state is $\rho _{1}$ and the second outcome never occurs. Thus
the average probability of error for this decision rule is 
\begin{equation}
R=\frac{1}{2}\left\langle \psi _{0}\right| \rho _{1}\left| \psi
_{0}\right\rangle ^{N}.
\end{equation}

Thus the rates of error decline coincide for the cases of joint and
separable measurements. Since the optimal separable test cannot do better
than the optimal joint measurement, the measurement considered is optimal.

QED.

If both states are mixed, then we can use the measurement that maximizes
fidelity distance between distributions of outcomes. In other words, the
measurement is chosen in such a way that it minimizes%
\begin{equation}
F(P,Q)=\sum \sqrt{p_{i}q_{i}}.
\end{equation}%
We will call this measurement fidelity-optimal. The advantage of this method
is that the fidelity-optimal measurement is easy to compute. It is simply a
measurement with outcomes that are orthogonal projectors on the eigenvectors
of the following operator:%
\begin{equation}
M=\rho _{1}^{-1/2}\sqrt{\rho _{1}^{1/2}\rho _{0}\rho _{1}^{1/2}}\rho
_{1}^{-1/2}.
\end{equation}%
(See \citeN{fuchs_caves95} for an explanation why this $M$ is
fidelity-optimal.) This measurement has only $d$ outcomes and their
probabilities are easy to compute. For this fidelity-optimal measurement we
can write a bound on the asymptotic error:

\begin{theorem}
The asymptotic error of the test based on the fidelity-optimal measurement
has the following asymptotic bound:%
\begin{equation}
\frac{1}{N}\log R\lesssim \log F(\rho _{0},\rho _{1}).
\end{equation}
\end{theorem}

This is the same upper bound that we have for joint asymptotic measurement
according to Theorem \ref{joint_bounds_theorem}.

\textbf{Proof:} 
\begin{equation}
\frac{1}{N}\ln R=\min_{0\leq \lambda \leq 1}\log
\sum_{i=1}^{N}p_{i}^{\lambda }q_{i}^{1-\lambda }\leq \log \sum_{i=1}^{N}%
\sqrt{p_{i}q_{i}}\leq \log F(\rho _{0},\rho _{1}).
\end{equation}%
The equality holds because of (\ref{classical_error_rate}), and the second
inequality is inequality (44) in \citeN{fuchs_graaf99}.

QED.

\section{Illustration}

This section illustrates the concepts developed above with an example of
testing for the presence of entanglement. Entanglement is one of the
properties of quantum systems that clearly separates them from classical
systems. It is a co-dependence of two remote parts of a quantum system that
cannot be created or destroyed by local operations on the parts.
Entanglement has become an important part of many quantum technologies
including quantum teleportation and quantum cryptography.

Entanglement has been produced in the laboratory. For example, %
\shortciteN{turchette98} developed a technique in which two ions are trapped
and illuminated equally by a laser beam that results in the creation of
entanglement.

In this illustration we are interested in tests of whether the entanglement
has been produced or not.

An example of an entangled quantum state is a pure state of the system of
two particles that corresponds to the projector on the following vector:%
\begin{equation}
\psi _{0}=\frac{1}{\sqrt{2}}\left( \left| 00\right\rangle +\left|
11\right\rangle \right) ,
\end{equation}%
where $\left| 00\right\rangle $ and $\left| 11\right\rangle $ denote $\left|
0\right\rangle \otimes \left| 0\right\rangle $ and $\left| 1\right\rangle
\otimes \left| 1\right\rangle ,$ and $\left| 0\right\rangle $ and $\left|
1\right\rangle $ form an orthonormal basis in the Hilbert space
corresponding to one of the particles.

The density matrix for this system is 
\begin{equation}
\rho _{0}=\left| \psi _{0}\right\rangle \left\langle \psi _{0}\right|
=\left( 
\begin{array}{cccc}
\frac{1}{2} & 0 & 0 & \frac{1}{2} \\ 
0 & 0 & 0 & 0 \\ 
0 & 0 & 0 & 0 \\ 
\frac{1}{2} & 0 & 0 & \frac{1}{2}%
\end{array}%
\right) .
\end{equation}%
The alternative hypothesis is that the state is a mix of two non-entangled
states given by projectors on vectors $\left| 00\right\rangle $ and $\left|
11\right\rangle ,$ respectively. The density matrix for this hypothesis is 
\begin{equation}
\rho _{1}=\left( 
\begin{array}{cccc}
\frac{1}{2} & 0 & 0 & 0 \\ 
0 & 0 & 0 & 0 \\ 
0 & 0 & 0 & 0 \\ 
0 & 0 & 0 & \frac{1}{2}%
\end{array}%
\right) .
\end{equation}%
This state can be easily produced by local operations but it is useless for
technologies that require entanglement.

Applying Theorem \ref{separable_pure}, we obtain the following formula for
the asymptotic error%
\begin{equation}
R\sim \frac{1}{2}\left\langle \psi _{0}\right| \rho _{1}\left| \psi
_{0}\right\rangle ^{N}=\frac{1}{2^{N+1}}.
\end{equation}%
It follows that it is sufficient to measure a sample of size $3$ to reduce
error below $5\%$.

The components of the optimal separable measurement are the projection on $%
\psi _{0}$ and its complement. Note that this is a joint measurement of both
particles. Actually, $\rho _{0}$ and $\rho _{1}$ cannot be distinguished by
the measurements that operate on each particle separately. This problem is
statistically unidentified by local measurements.

\section{Conclusion}

We have estimated the Chernoff efficiency bound for cases of joint and
separable measurements and also calculated it exactly for both pure and
commuting states. The results suggest that the loss of efficiency caused by
restriction to separable measurements is small.

Several questions remain open. Notably, it is not known whether the joint
measurement can ever be asymptotically better than the optimal separable
measurement. Second, it is not clear whether the optimal separable
measurement consists of orthogonal projectors. Third, it is unclear whether
the number of outcomes in this measurement is finite for finite-dimensional
quantum states.

\bibliographystyle{CHICAGO}
\bibliography{comtest}

\end{document}